\documentclass[twocolumn]{article}
\usepackage{authblk}
\usepackage{amsmath}
\usepackage{amssymb}
\usepackage{mathtools}
\usepackage[pdftex]{graphicx}
\usepackage[margin=20truemm]{geometry}

\usepackage{mathtools}
\usepackage{physics}
\usepackage{xcolor}
\usepackage{amsthm}
\newtheorem{defi}{Definition}
\newtheorem{prop}[defi]{Proposition}

\author[1]{Yasuhiro Matsumoto${}^{*}$}
\author[2]{Kei Matsushima}
\affil[1]{\small Center for Information Infrastructure, Institute of Science Tokyo}
\affil[2]{\small Institute of Engineering Innovation, Graduate School of Engineering, The University of Tokyo}
\affil[$*$]{\small Corresponding author}
\date{}
\title{Injectivity of boundary integral operator in direct-indirect mixed Burton-Miller equation for wave scattering problems with transmissive circular inclusion}

\begin{document}

\maketitle

\abstract{
This study proves that the injectivity condition for the integral operator of the direct-indirect
mixed Burton-Miller (BM) boundary integral equation (BIE) for Helmholtz transmission problems
is identical to that for the ordinary BM BIE for Helmholtz transmission problems
with a transmissive circular inclusion.
Although some numerical methods based on the direct-indirect mixed BM BIE
can be computed faster than the ordinary BM BIE, its well-posedness has been unclear.
This study resolves a part of the well-posedness, namely the injectivity of the integral operator with a transmissive circular inclusion.
}

\noindent {\bf Keywords}: Boundary integral equation, Burton-Miller method, transmission problem

\section{Introduction}
Wave scattering by transmissive scatterers is an important phenomenon in physics and engineering (e.g.\@, the design of metamaterials).
The simplest boundary value problem formulation of this phenomenon is the Helmholtz transmission problem \cite{kress1978}.
The conditions under which an analytical solution can be obtained for this problem are limited.
Therefore, a fast and accurate numerical method is required.
If the computational domain includes infinity,
the boundary integral equation (BIE) method, whose solution satisfies the radiation condition, is suitable.

At certain complex-valued frequencies, the Helmholtz transmission problem is known to have non-trivial (non-zero) solutions \cite{kress1978}.
These frequencies are called (true) eigenfrequencies.
The BIE method induces fictitious eigenfrequencies due to the formulation.
The Burton-Miller (BM) BIE \cite{burton1971application} ensures
the uniqueness of the solution at \textcolor{black}{real positive frequencies}.
Therefore, its extension to a transmission version is promising.
Moreover, the direct-indirect mixed BM BIE variant
shows better computational performance in \textcolor{black}{
  the fast direct solver based on the proxy-surface low-rank approximation method
  owing to the sparsity of the boundary integral operator configuration \cite{matsumoto2023jascome}.
}
However, the well-posedness of this variant is unclear.
A previous study \cite{rapun2008mixed} also used a direct-indirect mixed formulation to solve the transmission problem. However, it did not use the BM method; instead, it used an indirect integral representation to the exterior domain,
which suffers from fictitious eigenfrequencies at real frequencies.

The present study proves a part of the well-posedness of the mixed BM BIE.
Specifically, it proves that the injectivity condition of the integral operator of the mixed BM BIE is same as that of the BM BIE in case of a two-dimensional circular scatterer owing to Graf's addition theorem.

\section{Helmholtz transmission problem}
Let $\Omega_{1} \subset \mathbb R^2$ be a cylindrical scatterer with radius $a$.
The boundary of $\Omega_1$ is denoted by $\Gamma$.
We define the exterior infinite domain $\Omega_0 := \mathbb{R}^2 \setminus \overline{\Omega_{1}}$.
The positive constants $\varepsilon_{0}$, $\varepsilon_{1}$ and $\mu_{0}$, $\mu_{1}$ are material parameters for $\Omega_{0}$, $\Omega_{1}$, respectively.
The wavenumber in $\Omega_{j}$ is defined as $k_{j} := \omega \sqrt{\varepsilon_{j} \mu_{j}}$ for $j = 0, 1$, where $\omega$ is the angular frequency of incident wave $u^I$.
The wave scattering problem for transmissive scatterer $\Omega_1$ is to find the solution $u$ of a boundary value problem described as
\begin{align}
  &\Delta u(x) + k_{0}^{2} u(x) = 0, \quad x \in \Omega_{0}, \label{eq:bvp_start} \\
  &\Delta u(x) + k_{1}^{2} u(x) = 0, \quad x \in \Omega_{1}, \\
  &\lim_{h \downarrow 0} u(x + hn(x)) = \lim_{h \downarrow 0} u(x - hn(x)), \quad x \in \Gamma, \\
  &\frac{1}{\varepsilon_0} \lim_{h \downarrow 0} n(x) \cdot \nabla u(x + hn(x)) \nonumber \\
  &= \frac{1}{\varepsilon_1} \lim_{h \downarrow 0} n(x) \cdot \nabla u(x - hn(x)), \quad x \in \Gamma, \\
  &\lim_{|x| \to \infty} \sqrt{|x|} \qty(\pdv{}{|x|} - i k_{0}) (u - u^{I})(x) = 0, \label{eq:bvp_end}
\end{align}
where $n(x)$ is an outward unit normal vector at $x \in \Gamma$.

We define the following two potentials:
\begin{align}
  U_{0}(x) :=& u^{I}(x) - \varepsilon_{0} \int_{\Gamma} G^{k_0}(x - y) q(y) \dd{s(y)} \nonumber \\
  &+ \int_{\Gamma} \pdv{G^{k_0}}{n(y)} \qty(x - y) u(y) \dd{s(y)}, \quad x \notin \Gamma, \label{eq:bcbmex} \\
  U_{1}(x) :=& \varepsilon_{1} \int_{\Gamma} G^{k_1}(x - y) q(y) \dd{s(y)} \nonumber \\
  &- \int_{\Gamma} \pdv{G^{k_1}}{n(y)} \qty(x - y) u(y) \dd{s(y)}, \quad x \notin \Gamma, \label{eq:bcbmin}
\end{align}
where $G^{k}$ is the fundamental solution of the two-dimensional Helmholtz equation with wavenumber $k$, expressed as
\begin{align}
  G^{k}(x - y) := \frac{{i}}{4} H_{0}^{(1)} (k |x - y|), \quad x, y \in \mathbb{R}^2. \label{eq:funda}
\end{align}
Here, $i$ is the imaginary unit and $H_{m}^{(1)} (z)$ is the $m$-th order Hankel function of the first-kind.
In \eqref{eq:bcbmex} and \eqref{eq:bcbmin},
$\pdv{}{n(x)}$ is the normal derivative at $x$,
$q(x)$ is defined as
\begin{align}
  q(x) := \frac{1}{\varepsilon_0} \qty(\pdv{u}{n(x)})^{+} = \frac{1}{\varepsilon_1} \qty(\pdv{u}{n(x)})^{-},
\end{align}
where superscripts $+$ and $-$ are the limits from $\Omega_{0}$ and $\Omega_{1}$, respectively.

\section{Boundary integral equations}
\subsection{Burton-Miller boundary integral equation}
From Green's identities and the property of the fundamental solution, $U_{0}$ and $U_{1}$ respectively satisfy
\begin{align}
  U_{0}(x) = 
 \begin{cases}
    u(x), & x \in \Omega_0 \\
    0, & x \in \Omega_1
  \end{cases}
 , \quad
  U_{1}(x) = 
 \begin{cases}
    0, & x \in \Omega_0 \\
    u(x), & x \in \Omega_1
  \end{cases}
 .
\end{align}
\textcolor{black}{
In what follows, the incident wave term $u^I$ is neglected, as this paper focuses on non-trivial solutions corresponding to the eigenfrequencies.
Considering the limit of $U_{0}(x)$ to $\Gamma$ from $\Omega_{1}$, a linear combination of this limit and its normal derivative results in
\begin{align}
  &\qty( \qty(D_{}^{k_0} - \tfrac{1}{2}) + \alpha N_{}^{k_0}) u_{}(x) \notag \\
  &- \varepsilon_{0} \qty(S_{}^{k_0} + \alpha \qty(D_{}^{* k_0} + \tfrac{1}{2})) q_{}(x) = 0, \quad x \in \Gamma, \label{eq:BIEbmEx}
\end{align}
where $\alpha = i/k_{0}$ is a constant.
Similarly, the limit of $U_{1}(x)$ to $\Gamma$ from $\Omega_{0}$ results in
\begin{align}
  -\qty(D_{}^{k_1} + \tfrac{1}{2}) u_{}(x) + \varepsilon_{1} S_{}^{k_1} q_{}(x) = 0, \quad x \in \Gamma. \label{eq:BIEbmIn}
\end{align}
Coupling \eqref{eq:BIEbmEx} and \eqref{eq:BIEbmIn} yields the BM-formulated \cite{burton1971application} BIEs for \eqref{eq:bvp_start}--\eqref{eq:bvp_end}.
}
In \eqref{eq:BIEbmEx} and \eqref{eq:BIEbmIn},
$S^{k}$ and $D^{k}$ are the single- and double-layer potentials with wavenumber $k$, respectively.
$D^{\ast k}$ and $N^{k}$ are the normal derivatives of $S^{k}$ and $D^{k}$, respectively.
BIEs \eqref{eq:BIEbmEx} and \eqref{eq:BIEbmIn} can be rearranged into matrix form as
\begin{align}
  \mqty(
  \qty(D_{}^{k_0} - \tfrac{1}{2}) + \alpha N_{}^{k_0} & - \varepsilon_{0} \qty(S_{}^{k_0} + \alpha \qty(D_{}^{* k_0} + \tfrac{1}{2})) \\
  -\qty(D_{}^{k_1} + \tfrac{1}{2}) & \varepsilon_{1} S_{}^{k_1}
  ) \notag \\
  \times
  \mqty(
  u_{}(x) \\
  q_{}(x)
  )
  =
  \mqty(
  0 \\
  0
  ). \label{eq:BM_mat}
\end{align}

\subsection{Direct-indirect mixed Burton-Miller boundary integral equation}
Let $\phi$ be a density function on $\Gamma$ that satisfies the following indirect integral representation of a solution:
\begin{align}
  u(x) = \int_{\Gamma} G^{k_1} (x - y) \phi(y) \dd{s(y)}, \quad x \in \Omega_{1}.
\end{align}
Then, the BIE with respect to \eqref{eq:BIEbmIn} is replaced by the following indirect forms:
\begin{align}
  &u(x) = S^{k_1} \phi(x), &\quad x \in \Gamma, \label{eq:indirectBIE} \\
  &q(x) = \tfrac{1}{\varepsilon_{1}} \qty(D^{\ast k_1} + \tfrac{1}{2}) \phi(x), &\quad x \in \Gamma. \label{eq:inderectBIEdiff}
\end{align}
\textcolor{black}{Coupling \eqref{eq:indirectBIE}, \eqref{eq:inderectBIEdiff} and \eqref{eq:BIEbmEx}}
yields the direct-indirect mixed BM BIE \cite{matsumoto2023jascome} expressed as
\begin{align}
  \mqty(
  \begin{aligned}
    \qty(D^{k_0} - \tfrac{1}{2}) \\
    + \alpha N^{k_0}
  \end{aligned}
  & -\varepsilon_{0} \qty(S^{k_0} + \alpha \qty(D^{\ast k_0} + \tfrac{1}{2})) & 0 \\
  -1 & 0 & S^{k_1} \\
  0 & -1 & \tfrac{1}{\varepsilon_{1}} \qty(D^{\ast k_1} + \tfrac{1}{2})
  ) \nonumber \\
  \times
  \mqty(
  u(x) \\
  q(x) \\
  \phi(x)
  )
  =
  \mqty(
  0 \\
  0 \\
  0
  ), \quad x \in \Gamma. \label{eq:modiBM}
\end{align}

\section{Main result}
Before we state the proposition,
we define the eigenfrequencies of the integral operators.
\begin{defi}
  Let $A(\omega)$ be an integral operator such as the left-hand-side coefficient of \eqref{eq:BM_mat} or \eqref{eq:modiBM}.
  Then, the angular frequency $\omega$ is called the eigenfrequency or eigenvalue if $A(\omega)z=0$ has a non-trivial solution $z \neq 0$ at $\omega$.
\end{defi}
A part of the well-posedness of the BIE, in this case the injectivity of the integral operator, is related to the eigenfrequencies, as given in the following proposition.
\begin{prop} \label{prop:injectivity}
  Let $\Omega_{1}$ be a cylinder with radius $a$.
  Then, the condition for $\omega$ under which the integral operator of the direct-indirect mixed BM BIE \eqref{eq:modiBM} is injective is equivalent to the condition under which the integral operator of the BM BIE \eqref{eq:BM_mat} is injective.
  In other words, the distribution of eigenfrequencies of the direct-indirect mixed BM BIE is equivalent to that of the BM BIE.
\end{prop}
\begin{proof}
Since $u$ and $q$ are continuous on $\Gamma$, they can be Fourier-expanded to
\begin{align}
  u = \sum_{n = -\infty}^{\infty} u_{}^{n} e^{in\theta}, \\
  q = \sum_{n = -\infty}^{\infty} \textcolor{black}{q_{}^{n} e^{in\theta}},
\end{align}
in the circumferential direction, where $\theta$ is the polar coordinate angle and
$u^{n}$ and $q^{n}$ are the coefficients of the Fourier series expansion.
Applying Graf's addition theorem to the Hankel function in \eqref{eq:BIEbmEx}, as in \cite{misawa2012jascome} or as the
\textcolor{black}{one-sided} limit in \cite{klockner2013quadrature} results in
\begin{multline}
  \qty( \qty(D_{}^{k_0} - \tfrac{1}{2}) + \alpha N_{}^{k_0}) u_{}(x) - \varepsilon_{0} \qty(S_{}^{k_0} + \alpha \qty(D_{}^{* k_0} + \tfrac{1}{2})) q_{}(x) \\
  = \frac{i \pi a}{2} \sum_{n = -\infty}^{\infty} \left(k_{0} {H_{n}^{(1)}}^{\prime}(k_{0} a_{}) J_{n}(k_{0} a_{}) \right. \\
  + \left. \alpha k_{0}^{2} {H_{n}^{(1)}}^{\prime}(k_{0} a_{}) {J_{n}}^{\prime}(k_{0} a_{}) \right) u_{}^{n} e^{in\theta}  \\
  -\frac{i \pi a}{2} \sum_{n = -\infty}^{\infty} \varepsilon_{0} \left({H_{n}^{(1)}}^{}(k_{0} a_{}) J_{n}(k_{0} a_{}) \right. \\
  + \left. \alpha k_{0} {H_{n}^{(1)}}^{}(k_{0} a_{}) {J_{n}}^{\prime}(k_{0} a_{}) \right) q_{}^{n} e^{in\theta}
  = 0. \label{eq:wayEx}
\end{multline}
Applying the same calculation to \eqref{eq:BIEbmIn} yields
\begin{align}
  &-\qty(D_{}^{k_1} + \tfrac{1}{2}) u_{}(x) + \varepsilon_{} S_{}^{k_1} q_{}(x) = \nonumber \\
  &-\frac{i \pi a}{2} \sum_{n = -\infty}^{\infty} \qty(k_{1} H_{n}^{(1)}(k_{1} a_{}) {J_{n}}^{\prime}(k_{1} a_{})) u_{}^{n} e^{in\theta} \nonumber \\
  &+\frac{i \pi a}{2} \sum_{n = -\infty}^{\infty} \varepsilon_{1} \qty({H_{n}^{(1)}}^{}(k_{1} a_{}) J_{n}(k_{1} a_{})) q_{}^{n} e^{in\theta} = 0. \label{eq:wayIn}
\end{align}
Coupling \eqref{eq:wayEx} and \eqref{eq:wayIn}, we obtain the following matrix form:
\begin{align}
  \frac{i \pi a}{2} \sum_{n = -\infty}^{\infty}
  \mqty(
  J_{n}(k_{0} a_{}) + \alpha k_{0} {J_{n}}^{\prime}(k_{0} a_{}) & 0 \\
  0 & H_{n}^{(1)}(k_{1} a_{})
  ) \nonumber \\
  \times
  \mqty(
  k_{0} {H_{n}^{(1)}}^{\prime}(k_{0} a_{}) & -\varepsilon_{0} H_{n}^{(1)}(k_{0} a_{}) \\
  -k_{1}{J_{n}}^{\prime}(k_{1} a_{}) & \varepsilon_{1} J_{n}(k_{1} a_{})
  )
  \mqty(
  u_{}^{n} e^{in\theta} \\
  q_{}^{n} e^{in\theta}
  )
  =
  \mqty(
  0 \\
  0
  ). \label{eq:bm_fourier}
\end{align}

Since the set $\{e^{i n \theta} \}_{n = -\infty}^{\infty}$ is linearly independent, the existence of non-trivial solutions to the above equation must be considered independently for each $n$.
For the $n$-th order of the above equation,
the determinant on the left-hand-side vanishes
when either of the following conditions is satisfied:
\begin{align}
  -\varepsilon_{0} k_{1} H_{n}^{(1)}(k_{0} a_{}) {J_{n}}^{\prime}(k_{1} a_{}) + \varepsilon_{1} k_{0} {H_{n}^{(1)}}^{\prime}(k_{0} a_{}) {J_{n}}^{}(k_{1} a_{}) = 0, \label{eq:trueEigen}
\end{align}
\begin{align}
  H_{n}^{(1)}(k_{1} a_{}) \qty(J_{n}(k_{0} a_{}) + \alpha k_{0} {J_{n}}^{\prime}(k_{0} a_{})) = 0, \label{eq:fictitiousEigen}
\end{align}
where $k_{0} = \omega \sqrt{\varepsilon_{0} \mu_{0}}$ and $k_{1} = \omega \sqrt{\varepsilon_{1} \mu_{1}}$.

Similar to the scheme for the ordinary BM BIE,
we apply Fourier expansion and Graf's addition theorem to \eqref{eq:modiBM} to clarify the eigenvalues of \eqref{eq:modiBM}.
First, $\phi$ is Fourier-expanded to $\phi = \phi^{n} e^{in\theta}$.
Then, \eqref{eq:indirectBIE} is expanded to
\begin{align}
  \sum_{n = -\infty}^{\infty} u^{n} e^{i n \theta} = \frac{i \pi r}{2} \sum_{n = -\infty}^{\infty} H_{n}^{(1)}(k_{1}a) J_{n}(k_{1}a) \phi^{n} e^{i n \theta},
\end{align}
and \eqref{eq:inderectBIEdiff} is expanded to
\begin{align}
  \sum_{n = -\infty}^{\infty} q^{n} e^{i n \theta} = \frac{i \pi r}{2} \frac{k_{1}}{\varepsilon_{1}} \sum_{n = -\infty}^{\infty} H_{n}^{(1)}(k_{1}a) {J_{n}}^{\prime}(k_{1}a) \phi^{n} e^{i n \theta}.
\end{align}
Substituting the above two equations into \eqref{eq:modiBM}, \textcolor{black}{we have}
\begin{multline}
  \textcolor{black}{\frac{i \pi a}{2} \sum_{n = -\infty}^{\infty}} \\
  \mqty(
  \begin{aligned}
    &k_{0} {H_{n}^{(1)}}^{\prime}(k_{0}a) \\
    &\times \left\{ J_{n}(k_0 a) \right. \\
    &\left. + \alpha k_{0} {J_{n}}^{\prime}(k_0 a) \right\}
  \end{aligned}
  &
  \begin{aligned}
    &-\varepsilon_{0} {H_{n}^{(1)}}(k_{0}a) \\
    & \times \left\{ J_{n}(k_0 a) \right. \\
    & + \left. \alpha k_{0} {J_{n}}^{\prime}(k_0 a) \right\}
  \end{aligned}
  & 0 \\
  -1 & 0 & 
  \begin{aligned}
    H_{n}^{(1)}(k_{1}a) \\
    \times J_{n}(k_1 a)
  \end{aligned}
  \\
  0 & -1 &
  \begin{aligned}
  \tfrac{1}{\varepsilon_{1}} H_{n}^{(1)}(k_{1}a) \\
       \times {J_{n}}^{\prime}(k_1 a)
  \end{aligned}
  ) \\
  \times
  \mqty(
  u^{n} e^{i n \theta} \\
  q^{n} e^{i n \theta} \\
  \phi^{n} e^{i n \theta}
  )
  =
  \mqty(
  0 \\
  0 \\
  0
  ). \label{eq:modi_bm_fourier}
\end{multline}
\textcolor{black}{
  Using the orthogonality of $\{e^{in\theta}\}_{n = -\infty}^{\infty}$ and the rule of Sarrus, we have
}
\begin{multline}
  -\qty{-\varepsilon_{0} {H_{n}^{(1)}}^{\prime}(k_{0}a)
    \qty(J_{n}(k_0 a) + \alpha k_{0} {J_{n}}^{\prime}(k_0 a))} \\
  \times \qty(-\tfrac{1}{\varepsilon_{1}} H_{n}^{(1)}(k_{1}a) {J_{n}}^{\prime}(k_1 a)) \\
  - k_{0} {H_{n}^{(1)}}^{\prime}(k_{0}a) \qty( J_{n}(k_0 a)
  + \alpha k_{0} {J_{n}}^{\prime}(k_0 a)) \\
  \times \qty(-H_{n}^{(1)}(k_{1}a) J_{n}(k_1 a)) \\
  = \qty{H_{n}^{(1)}(k_{1}a) \qty(J_{n}(k_0 a) + \alpha k_{0} {J_{n}}^{\prime}(k_0 a))} \\
  \times \qty{-\tfrac{\varepsilon_{0}}{\varepsilon_1} k_{1} H_{n}^{(1)}(k_{0}a) {J_{n}}^{\prime}(k_{1}a)
  + k_{0} {H_{n}^{(1)}}^{\prime}(k_{0}a) J_{n}(k_{1}a)}.
\end{multline}
Considering the zeros of the right-hand side of the above equation with respect to $\omega$,
the first and second sets of curly brackets coincide with \eqref{eq:fictitiousEigen}
and \eqref{eq:trueEigen}, respectively.
Thus, the distribution of eigenfrequencies for the mixed BM BIE is identical to that for the ordinary BM BIE.
\end{proof}

In fact, the zeros with respect to $\omega$ corresponding to \eqref{eq:trueEigen} are the eigenvalues of the boundary value problem.
On the other hand, those corresponding to \eqref{eq:fictitiousEigen} are called fictitious eigenvalues because
they are not associated with the well-posedness of the original boundary value problem.
The eigenvalues of the BM BIE are identical to those of a single boundary integral equation obtained by Misawa {\it et al.} \cite{misawa2016jascome}.

\section{Numerical example}
\begin{figure}[tb]
  \begin{minipage}[t]{0.48\linewidth}
    \centering
    \includegraphics[width=0.85\linewidth]{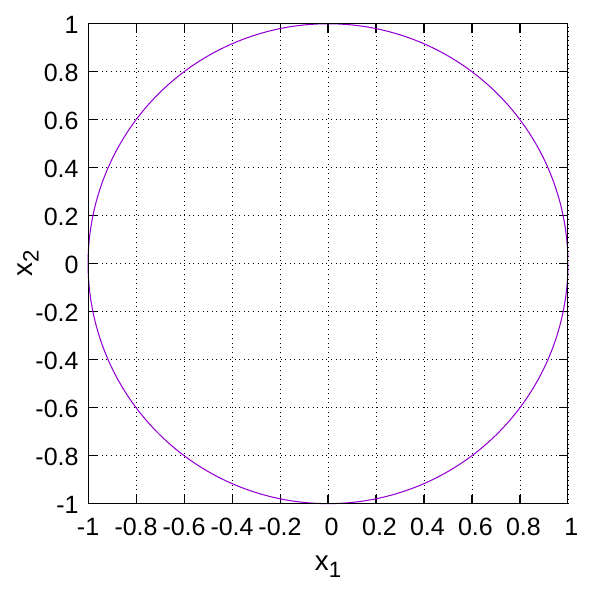}
    \caption{Circular shape $\Omega_{1}$.}
    \label{fig:circle}
  \end{minipage}
  \hfill
  \begin{minipage}[t]{0.48\linewidth}
    \centering
    \includegraphics[width=0.85\linewidth]{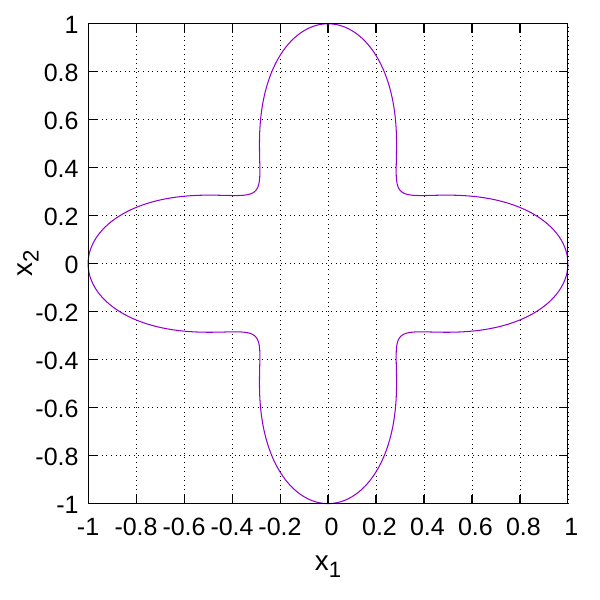}
    \caption{Star shape $\Omega_{1}$.}
    \label{fig:star}
  \end{minipage}
\end{figure}

We verify that Proposition \ref{prop:injectivity} holds numerically for the circular inclusion shown in Fig.\@ \ref{fig:circle}.
Then, we numerically show the possibility that Proposition \ref{prop:injectivity} holds even if the boundary shape of the inclusion is not circular (e.g.\@, the star shape shown in Fig.\@ \ref{fig:star}).

Since the coefficient matrix for the (discretized) BIE depends nonlinearly on the frequency $\omega$,
the Sakurai-Sugiura method \textcolor{black}{(SSM)} \cite{asakura2009} is employed to solve the nonlinear eigenvalue problem.
The BIEs of BM \eqref{eq:BM_mat} and direct-indirect mixed BM \eqref{eq:modiBM} are discretized using Nystr\"om methods based on \textcolor{black}{the zeta correction quadrature \cite{wu2023unified}}. %\cite{wu2021zeta, wu2023unified}.
In this quadrature, 31 local correction points are used, which results in an $O(N^{-31})$ convergence speed for $N$ global quadrature points.
In this section, material parameters $\varepsilon_{0} = 1$, $\varepsilon_{1} = 4$, and $\mu_{0} = \mu_{1} = 1$ and the number of global quadrature points $N = 400$ are fixed.
In the \textcolor{black}{SSM}, discretized linear equations for random right-hand-side vectors are solved using LU factorization.

\begin{figure}[tb]
  \centering
  \includegraphics[width=1.03\linewidth]{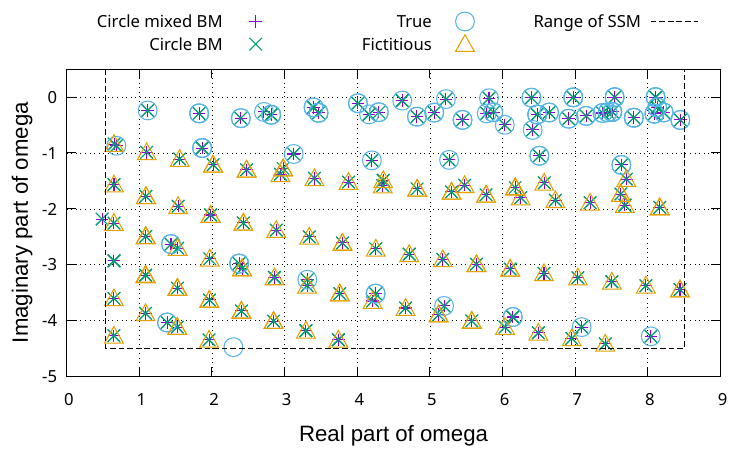}
  \caption{Eigenvalues computed using Sakurai-Sugiura method (SSM) for circular boundary.
    The eigenvalues of the mixed BM BIE well match those of the ordinary BM BIE.
    Furthermore, they are consistent with the zeros of either \eqref{eq:trueEigen} or \eqref{eq:fictitiousEigen},
    where ``True'' and ``Fictitious'' stand for the zeros of \eqref{eq:trueEigen} and \eqref{eq:fictitiousEigen}, respectively.
    \textcolor{black}{The range of SSM was divided into 4096 rectangles.
      On each side of a rectangle, the Gauss Legendre quadrature with 28 points were used.
      The number of SSM moments was 4.
      The eigenvalues near the real axis have non-zero imaginary part.}
  }
  \label{fig:ssm_circle}
\end{figure}

Fig.\@ \ref{fig:ssm_circle} shows plots of the complex-valued eigenfrequencies
computed by the \textcolor{black}{SSM} for a circular inclusion.
We can see that the eigenfrequencies of the mixed BM BIE well match those of the ordinary BM BIE.
Furthermore, these eigenfrequencies coincide with the zeros of either \eqref{eq:trueEigen} or \eqref{eq:fictitiousEigen}.
This result is consistent with the assertion of Proposition \ref{prop:injectivity}.

\begin{figure}[tb]
  \centering
  \includegraphics[width=1.03\linewidth]{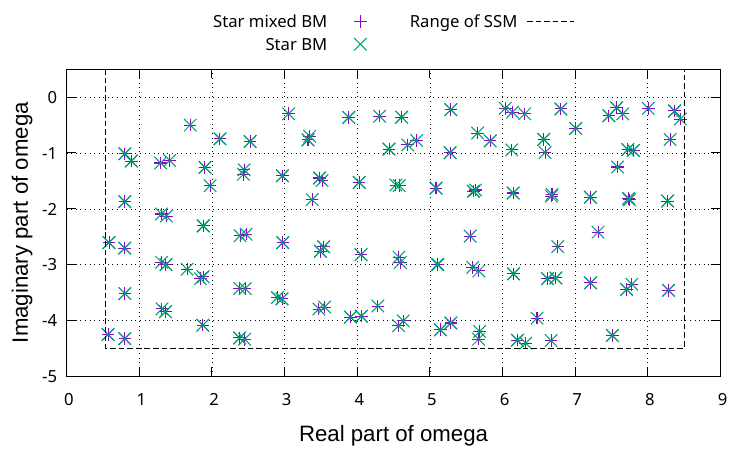}
  \caption{Eigenvalues computed using SSM for star-shaped boundary.
    The eigenvalues of the mixed BM BIE well match those of the ordinary BM BIE even though the inclusion is not circular.
    \textcolor{black}{The parameters of the SSM were the same as those of Fig. \ref{fig:ssm_circle}.}
  }
  \label{fig:ssm_star}
\end{figure}

Fig.\@ \ref{fig:ssm_star} shows plots of the \textcolor{black}{eigenvalues} %eigenfrequencies
computed by the \textcolor{black}{SSM} for a star-shaped inclusion.
Although the boundary shape is not a circle, it can be seen
that the eigenvalues of the mixed BM BIE are consistent with those of the ordinary BM BIE.
This result suggests the possibility that the assumptions of Proposition \ref{prop:injectivity} can be relaxed.

\section{Concluding remarks}
For Helmholtz transmission problems with a two-dimensional circular inclusion,
this paper proved that the condition of the \textcolor{black}{uniqueness for the solution of the direct-indirect mixed BM BIE}
coincides with that of the BM BIE
\textcolor{black}{if the solution exists.}
The correctness of the proposition was \textcolor{black}{also confirmed numerically in discretized system}.
\textcolor{black}{Numerical results} suggest that the proposition can be applied to \textcolor{black}{non-circular inclusions}.

\textcolor{black}{For establishing the well-posedness of the mixed BM BIE, we also require the surjectivity
on a suitable function space.%on a Banach space.
}
We will investigate the \textcolor{black}{bijectivity} under more general conditions.
\textcolor{black}{More numerical examples should be found where the mixed BM BIE shows excellent performance in numerical methods.}

\section*{Acknowledgments}
This work was supported by JSPS KAKENHI (grant numbers 23K19972, 24K20783 and 24K17191)
and computational resources of TSUBAME4.0 provided by Institute of Science Tokyo and Camphor3 provided by ACCMS of Kyoto University through JHPCN (jh250045).

\end{document}